\documentclass[11pt,reqno]{amsart}
\usepackage{mathrsfs}
\usepackage{amsmath}
\usepackage{amssymb, amsmath}
\usepackage{graphicx}

\newcommand{\newcom}{\newcommand}

\newcom{\al}{\alpha}
\newcom{\be}{\beta}
\newcom{\eps}{\epsilon}
\newcom{\veps}{\varepsilon}
\newcom{\ga}{\gamma}
\newcom{\Ga}{\Gamma}
\newcom{\ka}{\kappa}
\newcom{\Lam}{\Lambda}
\newcom{\lam}{\lambda}
\newcom{\vp}{\varphi}
\newcom{\Om}{\Omega}
\newcom{\om}{\omega}
\newcom{\Sig}{\Sigma}
\newcom{\sig}{\sigma}
\newcom{\tht}{\theta}
\newcom{\tri}{\triangle}
\newcom{\oo}{\infty}
\newcom{\vphi}{\varphi}
\newcom{\cB}{{\mathcal B}}
\newcom{\cC}{{\mathcal C}}
\newcom{\cD}{{\mathcal D}}
\newcom{\cF}{{\mathcal F}}
\newcom{\cL}{{\mathcal L}}
\newcom{\cM}{{\mathcal M}}
\newcom{\cP}{{\mathcal P}}
\newcom{\cS}{{\mathcal S}}
\newcom{\cQ}{{\mathcal Q}}
\newcom{\cT}{{\mathcal T}}
\newcom{\cY}{{\mathcal Y}}
\newcom{\cZ}{{\mathcal Z}}
\newcom{\R}{\Bbb R}
\newcom{\N}{\Bbb N}
\newcom{\Z}{\Bbb Z}
\newcom{\C}{\Bbb C}
\newcom{\E}{\Bbb E}
\def\dive{\mathop{\rm div}\nolimits}
\let\e=\varepsilon

\newcom{\f}{\frac}
\newcom{\dint}{\displaystyle\int}
\newcom{\dsum}{\displaystyle\sum}
\newcom{\dlim}{\displaystyle\lim}
\newcom{\ov}{\overline}
\newcom{\wt}{\widetilde}
\newcom{\pa}{\partial}
\newcom{\p}{\partial}
\newcom\na{\nabla}
\newcom\rto{\rightarrow}
\newcom\lto{\leftarrow}
\newcom\mto{\mapsto}
\newcom{\disp}{\displaystyle}
\newcom{\non}{\nonumber}
\newcom{\no}{\noindent}
\newcom{\QED}{$\square$}

\def\ef{\hphantom{MM}\hfill\llap{$\square$}\goodbreak}

\def\eqdefa{\buildrel\hbox{\footnotesize def}\over =}

\newtheorem{athm}{\bf \t}[section]
\newenvironment{thm} [1] {\def\t{#1}\begin{athm} \bf \rm} {\end{athm}}
\newcom{\bthm}{\begin{thm}}\newcom{\ethm}{\end{thm}}

\newcom{\beq}{\begin{equation}}
\newcom{\eeq}{\end{equation}}
\newcom{\ben}{\begin{eqnarray}}
\newcom{\een}{\end{eqnarray}}
\newcom{\beno}{\begin{eqnarray*}}
\newcom{\eeno}{\end{eqnarray*}}

\numberwithin{equation}{section}

\begin{document}

\title[Global Regularity for the
 Navier-Stokes equations] {Global Regularity for the Navier-Stokes equations with large, slowly varying initial data in the vertical direction}


\author{Marius Paicu}
\address{Departement de Mathematique, Universite Paris Sud, 91405 Orsay Cedex, FRANCE}
\email{Marius.Paicu@math.u-psud.fr}

\author{Zhifei Zhang}
\address{School of Mathematical Sciences, Peking University, 100871, P. R. China}
\email{zfzhang@math.pku.edu.cn}


\date{27,Mar}
%
%
%
%

\begin{abstract} In \cite{cgp} is obtained a class of large initial data generating a global smooth solution
to the three dimensional, incompressible Navier-Stokes equations.
This data varies slowly in the vertical direction (is a function on $\epsilon x_3$) and has a norm which blows up as the small parameter goes to zero. This type of initial data  can be seen as the ``ill prepared" case (in opposite with the ``well prepared" case which was treated  in  \cite{cgens}-\cite{cg3}). In \cite{cgp} the fluid evolves in a special domain, namely $\Omega=T^2_h\times\R_v$. The choice of a periodic domain in the horizontal variable plays an important role.
The aim of this article is to study the case where the fluid evolves in the full spaces $\R^3$, where we need to overcome the difficulties coming from  very low horizontal frequencies.  We consider in this paper an intermediate situation between the ``well prepared" case and ``ill prepared'' situation (the norms of the horizontal components of initial data are small but the norm of the vertical component  blows up as the small parameter goes to zero).
As in \cite{cgp}, the proof uses the analytical-type estimates and the special structure of the nonlinear term of the equation.
\end{abstract}


\maketitle


\setcounter{equation}{0}
\section{Introduction}\label{intro}
\setcounter{equation}{0}


We study in this paper the Navier-Stokes equations with initial data which is slowly varying in the vertical variable. More precisely we consider the system
\begin{equation*}(NS)
\begin{cases}
\partial_t u+u\cdot\nabla u- \Delta u=-\nabla p \quad\text{in}\quad \R_+\times \Omega \\
\dive u=0\\
u|_{t=0}=u_{0,\e},
\end{cases}
\end{equation*}
where $\Omega=\R^3$ and $u_{0,\e}$  is a divergence free vector field, whose dependence on the vertical variable~$x_3$ will be chosen to be ``slow'', meaning that it depends on~$\e x_3$ where~$\e$ is a small parameter.    Our goal is to prove a global existence in time result for the solution generated by this type of   initial data, with no smallness assumption on its norm. This type of initial data which is slowly varying in the vertical direction was already studied in \cite{cgens}-\cite{cg3}-\cite{cgp}. We recall that in \cite{cg3} is studied the case of ``well prepared" initial data of the type $(\e u^h_0(x_h,\e x_3), u^3_0(x_h,\e x_3))$ and in \cite{cgp} is studied the more difficult case of ``ill prepared" initial data of the type $(u^h(x_h,\e x_3),\e^{-1} u^3(x_h,\e x_3))$. In this paper, we consider the large initial data between the ``well prepared" case treated in \cite{cgens}-\cite{cg3} and the ``ill prepared" case treated in \cite{cgp}. More precisely, our initial data is of the form
$$u_{0,\e}=(\e^{\frac 12}u_0^h(x_h,\e x_3), \e^{-\frac 12} u_0^3(x_h,\e x_3)).$$

\smallskip

The mathematical study of the Navier-Stokes equations has a long history. We begin by recalling some important and classical  facts about the Navier-Stokes system, focusing on the conditions which imply  the global existence of the strong solution.

\smallskip

 The first important result about the classical Navier-Stokes system
was obtained by J. Leray \cite{leray}, and asserted that for  every finite energy initial data there exists at least one global in time weak solution which verifies the energy estimate. This solution is unique in $\R^2$ but unfortunately the solution is not known to be unique in three dimensional space. The result of J. Leray uses the structure of the nonlinear terms in order to obtain the energy inequality. The question of the uniqueness or of the regularity of the weak solutions is open.

The Fujita-Kato theorem gives a partial response
to the construction of global unique solution. Indeed, the theorem of Fujita-Kato \cite{fujitakato} allows to construct a unique  local in time solution in the homogeneous Sobolev spaces  $\dot H^{\frac 12}(\R^3)$, or in the Lebsegue space $L^3(\R^3)$
\cite{kato}. If the initial data is small compared to the viscosity  $\|u_0\|_{\dot H^{\frac 12}}\leq c\nu$, then the strong solution exists globally in time. This result was generalized by M. Cannone, Y. Meyer et  F. Planchon \cite{cannonemeyerplanchon} to Besov spaces of negative index of regularity. More precisely, they proved that, if the initial data belongs to the Besov space
$B^{-1+\frac 3p}_{p,\infty}(\R^3)$, and verifies  that it is small in the
norm of this Besov space, compared to the viscosity, then the
solution is global in time.  More recently, in \cite{kochtataru} is obtained a unique global in time solution for Navier-Stokes equation for small data belonging to a more general space of initial data, which is derivatives of $\text{BMO}$
function. Concerning the methods to obtain such results, we recall that proving the existence of a unique, global in time solution to the Navier-Stokes equations is
 rather standard  (it is a consequence of a Banach fixed point theorem) as long as the initial data is chosen small enough in some scale invariant spaces (with invariant  norm  by the scaling $\lambda u(\lambda^2t, \lambda x)$) embedded in~$\dot B^{-1}_{\infty,\infty}$ (the Besov space) , where we recall that
$$
\|f\|_{\dot B^{-1}_{\infty,\infty}} \eqdefa \sup_{t > 0} t^\frac12 \|e^{t\Delta} f\|_{L^\infty}.
$$
We refer for instance to~\cite{cannonemeyerplanchon},\cite{fujitakato},\cite{kochtataru},\cite{weisslerNS} for a proof in various scale invariant function spaces.  These theorems are general results of global existence for small initial data and does not take into account the any particular  algebraical properties of the nonlinear terms in the Navier-Stokes equations.

However, proving such a result without any smallness assumption or any  geometrical invariance hypothesis,
 which implies conservation of quantities beyond the scaling, is a challenge.  Little progress has been made
 in that direction: we will not describe all the literature on the question, but refer among others
 to~\cite{babinandco}, \cite{cannonemeyerplanchon},  \cite{cdgg2},  \cite{cdgg3}, \cite{cgens}-\cite{cg3}, \cite{cgp}, \cite{cz1}, \cite{raugelsell} and the references
therein for more details.

We recall briefly  the examples of large initial data which gives global
existence of the solution known in the literature. We first notice that for regular  axi-symmetric initial data,  without swirl, there exists a unique global in time solution for the Navier-Stokes system. This result is based on  the conservation of some quantities beyond the scaling regularity level (see \cite{ukhoiudo}).

The case of large initial data (in some sense) for  fluids evolving in thin domains was firstly considered by \cite{raugelsell}.  Roughly speaking, the tridimensional  Navier-Stokes system can be seen as a perturbation of the two dimensional Navier-Stokes system if the domain is thin enough in the vertical direction. Generally, if the initial data $u_0$ can be splitting as $u_0=v_0+w_0$, with $v_0$ a two dimensional free divergence vector field belonging to $L^2({\mathcal T}^2_h)$ and $w_0\in H^{\frac 12}(\mathcal T^3)$, such that
$$\|w_0\|_{H^{\frac 12}(\mathcal T^3)}\exp \bigg(\frac{\|v_0\|_{L^2(\mathcal T^2_h)}^2}{\nu^2}\bigg)\leq c\nu,$$ then the solution exists globally in time.

The case of initial data with  large initial vortex in the vertical direction ($\text{rot\,} u_0^\e=\text{rot\,}u_0+\epsilon^{-1} (0,0,1)$), or equivalently the case of rotating fluids,  was studied in \cite{babinandco} in the case of periodic domains and in \cite{cdgg2}-\cite{cdgg3} for the case of a rotating fluid in $\R^3$ or in $\R^2\times (0,1)$.  When the rotation is fast enough the fluid tends to have a two-dimensional behavior,
far from the boundary of the domain (this is the so called Taylor-Proudman
column theorem \cite{pedlovsky}). For example,  in the case where the domain is $\R^3$ the fluctuation of this motion is dispersed to infinity and some Strichartz quantities became small \cite{cdgg2} which allow to obtain the global existence of the solution (when $\e$ small enough).

An important issue for the Navier-Stokes equations is to use on maximum the algebrical structure of the nonlinear terms.  Some results used in a crucial way this structure and allow to obtain very interesting new results.

The case of the Navier-Stokes equations
with vanishing vertical viscosity was firstly studied in \cite{cdgg2} where the authors proved local existence for large data in anisotropic Sobolev spaces $H^{0,s}, s > 1/2$,
and global existence and uniqueness for small initial data. One
of the key observations is that, even if there is no vertical
viscosity and thus no smoothing in the vertical variable,
the partial derivative $\partial_3$ is only applied
to the component $u_3$ in the nonlinear term. The divergence-free
condition implies that $\partial_3 u_3$ is regular enough to get good
estimates of the nonlinear term. In \cite{cz1} the authors obtained the global existence of the solution for  the anisotropic Navier-Stokes system with high oscillatory initial data.

A different idea, but always using the special structure of the Navier-Stokes equations, is used by J.-Y. Chemin and I. Gallagher \cite{cgens} in order to construct the first example of periodic initial data which is big in $C^{-1}$, and strongly oscillating in one
direction which generates a global solution. Such initial data is given by
$$u_0^N=(Nu_h(x_h)\cos(Nx_3),-\dive_h u_h(x_h)\sin(Nx_3)),$$
where $\|u_h\|_{L^2(\mathcal T^2_h)}\leq C(\ln N)^{\frac 19}$. This result was generalized to the case of the space $\R^3$ in \cite{cg2}.

In the paper \cite{cg3}, J.-Y. Chemin and I. Gallagher  studied the Navier-Stokes equations for initial data which  slowly varies in the vertical direction in the well prepared case. The ``well prepared" case means that the norm of the initial data is large but  does not blow up when the parameter $\e$ converges to zero.  We note that  important remarks on the pressure term and the bilinear term were used in this paper in order to obtain the global existence for large data.

The case of slowly varying initial data in the vertical direction in the ``ill prepared" initial data was recently studied in \cite{cgp}.  We note that the horizontal components has a large norm and the vertical component has a norm which blows up when the parameter goes to zero. After a change of scale  of the problem, the system became a Navier-Stokes type equation with an anisotropic viscosity $-\nu\Delta_h u -\nu\epsilon^2\partial_3^2 u$ and anisotropic gradient of the pressure, namely $-(\nabla_h p, \epsilon^2 \partial_3 p)$. In this equation we can remark that there is a loos of regularity in the vertical variable in Sobolev estimates. To overcome this difficulty is needed to work with analytical initial data.
The most important tool was developed in the paper of J.-Y. Chemin \cite{chemin21} and consisted to make analytical type estimates, and in the same time to control the size of the analyticity band.  This is performed by the control of nonlinear quantities which depend on the solution itself. Even in this situation, it is important to take into account very carefully the special structure of the Navier-Stokes equations. In \cite{cgp} is obtained in fact a global in time Cauchy-Kowalewskaya type theorem. We recall also that  some local in time results for Euler and Prandtl equation with analytic initial data can be found in\cite{sammartino&caflisch}.

In \cite{cgp} the fluid is supposed to evolve in a special domain $\Omega=T^2_h\times\R_v$. This choice of domain is justified by the pressure term. Indeed, the pressure verifies the elliptic equation $\Delta_\e p=\partial_i\partial_j (u^i u^j)$, and consequently, $\nabla_h p=(-\Delta_\e)^{-1}\nabla_h \partial_i\partial_j(u^i u^j)$. Because we have that $\Delta_\e^{-1}$ converges to $\Delta_h^{-1}$ it is important to control the low horizontal frequencies. While in the case of the periodic torus in the horizontal variable we have only zero horizontal frequency and high horizontal frequencies.

In this paper our goal is to investigate the case where the fluid evolves in the full space $\R^3$. In that situation, we are able to solve globally in time the equation (conveniently rescaled in~$\e$) for small  analytic-type initial data. In the case of the full space $\R^3$ we need to control very precisely the low horizontal frequencies. We also note that we can construct functional  spaces where the operator $\Delta_h^{-1}\nabla_h(a\nabla_h b)$ is a bounded operator. However we still need to impose on the initial data more control of the regularity in the low horizontal frequencies (namely we impose that $u_0(\cdot, x_3)\in L^2(\R^2_h)\cap \dot H^{-\frac 12}(\R^2_h)$). In the vertical variable we need to impose analyticity of the data. The method of the proof follows closely the argument of \cite{cgp}, but instead to use pointwise estimates on the fourier variables, we write an equation with a regularizing term in the vertical variable and we use energy estimates on anisotropic Sobolev spaces of the form $H^{0,s}$ respectively $H^{-\frac12,s}$.

 Our main result in the case of the full space $\R^3$ is the following (for the notations see the next section).
\bthm{Theorem}  \label{thm-NS}
Let $a$ be a positive number, $s>\f12$. There exist two positive constants $\veps_0$ and $\eta$
such that for any divergence free fields $v_0$ satisfying
\beno
\|e^{a|D_3|}v_0\|_{H^{0,s}}+\|e^{a|D_3|}v_0\|_{H^{-\f12,s}}\le \eta,
\eeno
and for any $\veps\in (0,\veps_0)$, the Navier-Stokes system $(NS)$ with initial data
$$u_0^\e=\big(\e^{\frac 12}v^h_0(x_h,\e x_3),\e^{-\frac 12} v_0^3(x_h,\e x_3)\big)$$
has a global smooth solution on $\R^3$.
\ethm

As we already explain above, in order to prove the main theorem \ref{thm-NS}, we will firstly transform the system using the change of scale $$u^\e(t,x_h,x_3)=\big(\e^{\frac 12}v^h(t,x_h,\e x_3),\e^{-\frac 12}v^3(t,x_h,\e x_3)\big)$$
into a system of Navier-Stokes type, with a vertical vanishing viscosity, that is the Laplacian operator became $-\nu\Delta_h v-\epsilon^2\partial_3 v$ and a changed pressure term became $-(\nabla_h p,\epsilon^2\partial_3 p)$.

\medskip

Taking the advantage that we work in the full spaces $\R^3$, we can also consider a different type of initial data, with larger amplitude but strongly oscillating in the horizontal variables, namely, initial data of the form
$$u_0^\e=\big(\e^{-\frac 12}v_0^h(\e^{-1}x_h,x_3),\e^{-\frac 32}v_0^3(\e^{-1}x_h,x_3)\big).$$
However, this type of initial data has the $\dot B^{-1}_{\infty,\infty}$ norm on the same order as the initial data in the previous theorem. In order to solve the Navier-Stokes equations with this new type of initial data, we make the different change of scale,
$$
u^\e(t,x_h,x_3)=\big(\epsilon^{-\frac 12}v^h(\epsilon^{-2}t,\epsilon^{-1}x_h,x_3),\epsilon^{-\frac 32}v^3(\epsilon^{-2}t,\epsilon^{-1}x_h,x_3)\big)
$$
and we note that the rescaled system that we obtain is exactly the same as that in the proof of theorem \ref{thm-NS}. Consequently, we also obtain the following result.
\bthm{Theorem}\label{thm2-NS}
Let $a$ be a positive number, $s>\f12$. There exist two positive constants $\veps_0$ and $\eta$
such that for any divergence free fields $v_0$ satisfying
\beno
\|e^{a|D_3|}v_0\|_{H^{0,s}}+\|e^{a|D_3|}v_0\|_{H^{-\f12,s}}\le \eta,
\eeno
and for any $\veps\in (0,\veps_0)$, the Navier-Stokes system $(NS)$ with initial data
$$u_0^\e=\big(\e^{-\frac12}v^h_0(\e^{-1}x_h, x_3),\e^{-\frac 32} v_0^3(\e^{-1}x_h, x_3)\big)$$
has a global smooth solution on $\R^3$.
\ethm

\section{A simplified model}

 Let us consider the following equation
$$
\partial_t u+\gamma u+a(D) Q(u,u)=0$$
where $a(D)$ is a fourier multiplier of order one and, $Q$ is any quadratic form.
Then, if the initial data verifies
$$\|u_0\|_{X}=\int e^{\delta |\xi|} \hat u(\xi) d\xi\leq c\gamma\quad\text{whith}\quad a>0$$
then we have a global solution in the same space. We follow the method introduced in \cite{cgp} and \cite{chemin21}.
The idea of the proof  is the following, we want to control the some kind of quantities on the solution, but we must prevent the possible loose of the rayon of the analyticity of the solution. Let us introduce $\theta(t)$ the ``loose of analyticity", such that $\dot\theta (t)=\int e^{(\delta-\theta(t))|\xi|}|\hat u(\xi)|d\xi$, $\theta(0)=0$. We denote by $\Phi=(a-\lambda\theta(t))|\xi|$ and we define
$$\dot\theta(t)=\int |\hat u_\Phi(\xi)|d\xi=\|u_\Phi\|_X, \quad\theta(0)=0.$$
The computations which follow are performed under the condition   $\theta(t)\leq a/\lambda$ (which implies $\Phi\geq 0$).
The equation verified by  $\hat u_\phi$ is the following
$$
\partial_t \hat u_\Phi +\gamma \hat u_\Phi +\lambda \dot \theta(t)|\xi|\hat u_\Phi+a(\xi)e^{\Phi}(\widehat{u^2})=0.
$$
As $\dot\theta\geq 0$, after integration in $\xi$, we obtained the following equation
$$
\partial_t \|u_\Phi\|_X+\gamma \|u_\Phi\|_{X}+\lambda\dot\theta(t)\int |\xi||\hat u_\Phi| d\xi\leq C \int |\xi|  |\hat {u_\Phi}|\star |\hat u_\Phi| (\xi) d\xi.
$$
As $|\xi|\leq |\xi-\eta|+|\eta|$, we obtain
$$
\int |\xi|  \hat {u_\Phi}\star \hat u_\Phi (\xi) d\xi\leq 2\bigg(\int |\xi||\hat u_\phi|d\xi\bigg)\bigg(\int |\hat u_\Phi|d\xi\bigg)=2\dot\theta (t) \bigg(\int |\xi||\hat u_\phi|d\xi\bigg).
$$
So, choosing  $\lambda=4C$ we obtain
$$\dot\theta(t)=\|u_{\Phi}(t)\|_{X}\leq 2\|e^{a|D|} u_0\|_{X} e^{-\gamma t}$$
which, for  $u_0$ small enough, gives
$$
\theta(t)\leq \gamma^{-1}\|e^{a|D|}u_0\|_{X}\leq a\lambda^{-1}.
$$
This allows to obtain the global in time existence of the solution.



\section{Structure of the proof}

\subsection{Reduction to a rescaled problem}

We seek the solution of the form
\beno
u_{\veps}(t,x)\eqdefa \Bigl(\veps^\f12v^h(t,x_h,\veps x_3),{\veps^{-\f12}}v^3(t,x_h,\veps x_3)\Bigr).
\eeno
This leads to the following rescaled Navier-Stokes system
\begin{equation}\label{equ:RNS}
(RNS_\veps)\quad\left\{
\begin{array}{ll}
\p_tv^h-\Delta_h v^h-\veps^2 \p_3^2 v^h+\veps^\f12 v\cdot \na v^h=-\na^h q,\\
\p_tv^3-\Delta_h v^3-\veps^2 \p_3^2 v^3+\veps^\f12v\cdot\na v^3=-\veps^2\p_3 q,\\
\textrm{div}v= 0,\\
v(0)=v_{0}(x),
\end{array}
\right.
\end{equation}
where $\Delta_h\eqdefa \p_1^2+\p_2^2$ and $\na_h\eqdefa(\p_1,\p_2)$.
As there is no boundary, the rescaled pressure $q$ can be computed
with the formula
\ben\label{equ:pressure}
-\Delta_{\veps}q=\veps^\f12\textrm{div}_h(v\cdot\na v), \quad \Delta_\veps=\Delta_h+\veps^2\p_3^2.
\een

When $\veps$ tends to zero, $\Delta_\veps^{-1}$ looks like  $\Delta_h^{-1}$. Thus for low horizontal frequencies, an expression of $\na_h\Delta^{-1}_h$
cannot be estimated in $L^2$. This is one reason why the authors in \cite{cgp} work in $\Bbb{T}^2\times \R$.
To obtain a similar result in $\R^3$, we need to introduce the following anisotropic Sobolev space.

\bthm{Definition} Let $s,\sigma\in \R$, $\sigma<1$. The anisotropic Sobolev space $H^{\sigma,s}$ is
defined by
$$H^{\sigma,s}=\{f\in {\cS}'(\R^3); \|f\|_{H^{\sigma,s}}<\infty\},$$ where
$$
\|f\|_{H^{\sigma,s}}^2\eqdefa
\int_{\R^3}|\xi_h|^{2\sigma}(1+|\xi_3|^2)^s|\hat f(\xi)|^2d\xi,\quad
\xi=(\xi_h,\xi_3).
$$
For any $f,g\in H^{\sigma,s}$, we denote
\beno
(f,g)_{H^{\sigma,s}}\eqdefa (|D_h|^\sigma\langle D_3\rangle^sf,|D_h|^\sigma\langle D_3\rangle^sg)_{L^2},
\quad \langle D_3\rangle=(1+|D_3|^2)^\f12.
\eeno
\ethm

We prove that

\bthm{Theorem}\label{thm-RANS}
Let $a$ be a positive number, $s>\f12$. There exist two positive constants $\veps_0$ and $\eta$
such that for any divergence free fields $v_0$ satisfying
\beno
\|e^{a|D_3|}v_0\|_{H^{0,s}}+\|e^{a|D_3|}v_0\|_{H^{-\f12,s}}\le \eta,
\eeno
and for any $\veps\in (0,\veps_0)$, $(RNS_\veps)$ has a global smooth solution on $\R^3$.
\ethm

\subsection{Definition of the functional setting}

As in \cite{cgp}, the proof relies on exponential decay estimates
for the Fourier transform of the solution. Thus, for any locally bounded function $\Psi$ on $\R^+\times \R^3$ and for any
function $f$, continuous in time and compactly supported in Fourier space, we define
\beno
f_\Psi(t)\eqdefa \cF^{-1}\bigl(e^{\Psi(t,\cdot)}\widehat f(t,\cdot)\bigr).
\eeno
Now we introduce two key quantities we want to control in order to prove the theorem. We define the function $\theta(t)$ by
\ben\label{equ:theta}
\dot\theta(t)\eqdefa\veps\|v_\Phi^h(t)\|_{H^{\f12,s}}^2+\|v_\Phi^3(t)\|_{H^{\f12,s}}^2\quad
\textrm{and} \quad \theta(0)=0, \een
and we also define
\ben\label{equ:Psi}
\Psi(t)\eqdefa\|v_\Phi(t)\|^2_{H^{0,s}}+\int_0^t\|\nabla_hv_\Phi(\tau)\|^2_{H^{0,s}}d\tau,
\een
where
\ben\label{equ:Phi-def}
\Phi(t,\xi)\eqdefa (a-\lambda\theta(t))|\xi_3| \een for some
$\lambda$ that will be chosen later on.

\subsection{Main steps of the proof}

\bthm{Proposition}\label{prop:tht}
A constant $C_0$ exists such that, for any positive $\lambda$ and
for any $t$ satisfying $\theta(t)\le a/\lambda$, we have
\beno
\tht(t)\le\exp\bigl(C_0\Psi(t)\bigr)\Bigl[\|e^{a|D_3|}v_0\|_{H^{-\f12,s}}^2+
C_0\int_0^t\dot\tht(\tau)\Psi(\tau)d\tau\Bigr].
\eeno
\ethm

\bthm{Proposition}\label{prop:Psi}
There exist $C_1$ and $\lambda_0$ such that for $\lambda\ge \lambda_0$ and
for any $t$ satisfying $\theta(t)\le a/\lambda$, we have
\beno
\Psi(t)\le \|e^{a|D_3|}v_0\|_{H^{0,s}}^2\exp\bigl(C_1\Psi(t)\bigr).
\eeno
\ethm

The proof of Proposition \ref{prop:tht} and \ref{prop:Psi} will be presented in section 4 and section 5 respectively. For the moment,
let us assume that they are true and conclude the proof of Theorem \ref{thm-RANS}. As in \cite{cgp}, we use a continuation argument.
For any $\lambda\ge \lambda_0$ and $\eta$, let us define
\beno
\cT_\lambda\eqdefa\{T: \theta(T)\le 4\eta^2, \Psi(T)\le 2\eta^2\}.
\eeno
Similar to the argument in \cite{cgp}, $\cT_\lambda$ is of the form $[0,T^*)$ for some positive $T^*$.
Thus, it suffices to prove that $T^*=+\oo$. In order to use Proposition \ref{prop:tht} and \ref{prop:Psi}, we need to assume that
$\theta(T)\le \f a\lambda$, which leads to the condition
\beno
4\eta^2\le \f a\lambda.
\eeno
From Proposition \ref{prop:tht} and \ref{prop:Psi}, it follows that for all $T\in \cT_\lambda$,
\begin{equation} \label{equ:thtpsi}
\begin{split}
\theta(T)&\le\exp(2C_0\eta^2)(\eta^2+2C_0\eta^2\theta(T)),\\
\Psi(T)&\le \eta^2\exp(2C_1\eta^2)
 \end{split}
 \end{equation}
Now we choose $\eta$ such that
$$
\exp(2C_0\eta^2)<2, \quad \exp(2C_1\eta^2)<2,\quad 4C_0\eta^2<\f12.
$$
With this choice of $\eta$, then we infer from (\ref{equ:thtpsi})  that
\ben\label{thetasmall}
\theta(T)< 4\eta^2,\quad \Psi(T)<2\eta^2,
\een
which ensures that $T^*=+\oo$,
thus we conclude the proof of Theorem \ref{thm-RANS}.\ef

\section{The action of subadditive phases on products}

For any function $f$, we denote by $f^+$ the inverse Fourier
transform of $|\widehat f|$. Let us notice that the map $f\mapsto
f^+$ preserves the norm of all $H^{\sigma,s}$ spaces. Throughout
this section, $\Psi$ will denote a locally bounded function on
$\R^+\times \R^3$ which satisfies the following inequality
\ben\label{subaddi} \Psi(t,\xi)\le \Psi(t,\xi-\eta)+\Psi(t,\eta).
\een

Before presenting the product estimates, let us recall the Littlewood-Paley decomposition. Choose two
nonnegative even functions $\chi$, $\varphi \in {\cS}(\R)$ supported
respectively in ${\cB}=\{\xi\in\R,\, |\xi|\le\frac{4}{3}\}$ and
${\cC}=\{\xi\in\R,\, \frac{3}{4}\le|\xi|\le\frac{8}{3}\}$ such that
\beno
&&\chi(\xi)+\sum_{j\ge0}\varphi(2^{-j}\xi)=1,\quad\mbox{for}\quad \xi\in\R,\\
&&\sum_{j\in \Z}\varphi(2^{-j}\xi)=1, \quad\mbox{for}\quad
\xi\in\R\setminus\{0\}. \eeno The frequency localization operators
$\Delta_j^v$ and $S_j^v$ in the vertical direction are defined by
\beno
&&\Delta_j^vf=\cF^{-1}\bigl(\varphi(2^{-j}|\xi_3|)\widehat{f}\bigr)\quad\mbox{for}\,
j\ge 0,\quad
S_j^vf=\cF^{-1}\bigl(\chi(2^{-j}|\xi_3|)\widehat{f}\bigr)=\sum_{j'\le j-1}\Delta_{j'}^vf,\\
&&\Delta_{-1}^vf=S_{0}^vf, \quad\Delta_{j}^vf=0\quad\mbox{for}\,
j\le -2. \eeno And the frequency localization operators
$\dot\Delta_j^h$ and $S_j^h$ in the horizontal direction are defined
by \beno
&&\dot\Delta_j^hf=\cF^{-1}\bigl(\varphi(2^{-j}|\xi_h|)\widehat
f\bigr), \quad S_j^hf=\sum_{j'\le
j-1}\dot\Delta_{j'}^hf,\quad\mbox{for}\,\, j\in \Z. \eeno It is easy
to verify that
\ben\label{equ:equinorm}
\|f\|_{H^{\sigma,s}}^2\approx
\sum_{j,k\in\Z}2^{2js}2^{2k\sigma}\|\Delta_j^v\dot\Delta_k^hf\|_{L^2}^2.
\een

In the sequel, we will constantly use the Bony's decomposition from
\cite{Bony} that
\beq\label{VBony}
fg=T_f^vg+R_f^vg,
\eeq
with $$T_f^vg=\sum_{j}S_{j-1}^vf\Delta_j^vg, \quad R_fg=\sum_{j}S_{j+2}^vf \Delta_{j}^vg.$$
We also use the Bony's decomposition in the horizontal direction
\beq\label{HBony}
fg=T^h_fg+T_f^hg+R^h(f,g),
\eeq
with $$T^h_fg=\sum_{j}S_{j-1}^hf\dot\Delta_j^hg, \quad R^h(f,g)=\sum_{|j'-j|\le 1}\dot\Delta_j^hf\dot\Delta_{j'}^hg.$$

\bthm{Lemma}(Bernstein's inequality)\label{Lem:Berstein} Let $1\le
p\le q\le\infty$. Assume that $f\in L^p(\R^d)$, then there exists a
constant $C$ independent of $f$, $j$ such that \beno &&{\rm
supp}\hat{f}\subset\{|\xi|\leq C2^j\}\Rightarrow \|\partial^\alpha
f\|_{L^q}\le
C2^{j{|\alpha|}+dj(\frac{1}{p}-\frac{1}{q})}\|f\|_{L^p},
\\
&&{\rm supp}\hat{f}\subset\{\f1{C}2^j\leq |\xi|\leq C
2^j\}\Rightarrow \|f\|_{L^p}\le
C2^{-j|\alpha|}\sup_{|\beta|=|\al|}\|\partial^\beta f\|_{L^p}. \eeno
\ethm

\bthm{Lemma}\label{Lem:localproduct} Let
$s>\f12,\sigma_1,\sigma_2<1$ and $\sigma_1+\sigma_2>0$. Assume that
$a_\Psi\in H^{\sigma_1,s}$ and $b_\Psi\in H^{\sigma_2,s}$. Then
there holds \beno
\|\bigl[\Delta_j^v\dot\Delta_k^h(T_a^vb)\bigr]_\Psi\|_{L^2}+\|\bigl[\Delta_j^v\dot\Delta_k^h(R_a^vb)\bigr]_\Psi\|_{L^2}\le
Cc_{j,k}2^{(1-\sigma_1-\sigma_2)k}2^{-js}\|a_\Psi\|_{H^{\sigma_1,s}}\|b_\Psi\|_{H^{\sigma_2,s}},
\eeno with the sequence $(c_{j,k})_{j,k\in\Z}$ satisfying
$\disp\sum_{j,k}c_{j,k}\le 1$. \ethm

\no{\bf Proof.}\,\,Let us firstly prove  the case when the function
$\Psi$ is identically 0. Below we only present the proof of $R_ab$,
the proof for $T_ab$ is very similar. Using Bony's decomposition
(\ref{HBony}) in the horizontal direction, we write \beno
\Delta_j\dot\Delta_k^h(R_a^vb)&=&\sum_{j'}\Delta_j^v\dot\Delta_k^h(S_{j'+2}^va\Delta_{j'}^vb)\\
&=&\sum_{j'}\Delta_j^v\dot\Delta_k^h\bigl(T_{S_{j'+2}^va}^h\Delta_{j'}^vb
+T_{\Delta_{j'}^vb}^hS_{j'+2}^va+R^h(S_{j'+2}^va,\Delta_{j'}^vb)\bigr)\\
&:=& I+II+III. \eeno Considering the support of the Fourier
transform of $T_{S_{j'+2}^va}^h\Delta_{j'}^vb$, we have \beno
I=\sum_{j'\ge j-4}\sum_{|k'-k|\le 4}\Delta_j^v\dot\Delta_k^h
\bigl(S_{j'+2}^vS_{k'-1}^ha\Delta_{j'}^v\dot\Delta_{k'}^hb\bigr).
\eeno Then we get by Lemma \ref{Lem:Berstein} that \beno
\|I\|_{L^2}&\le&C\sum_{j'\ge j-4}\sum_{|k'-k|\le 4}\|S_{j'+2}^vS_{k'-1}^ha\Delta_{j'}\dot\Delta_{k'}^hb\|_{L^2}\\
&\le& C\sum_{j'\ge j-4}\sum_{|k'-k|\le
4}\|S_{j'+2}^vS_{k'-1}^ha\|_{L^\oo}\|\Delta_{j'}^v\dot\Delta_{k'}^hb\|_{L^2}.
\eeno We use Lemma \ref{Lem:Berstein} again to get \beno
\|S_{j'+2}^vS_{k'-1}^ha\|_{L^\infty}&\le& \sum_{j''\le
j'+1}\sum_{k''\le k'-2}
\|\Delta_{j''}^v\dot\Delta_{k''}^ha\|_{L^\oo}\\
&\le& C\sum_{j''\le j'+1}\sum_{k''\le k'-2}2^{k''}\|\Delta_{j''}^v\dot\Delta_{k''}^ha\|_{L^\oo_{x_3} L_{x_h}^2}\\
&\le& C\sum_{j''\le j'+1}\sum_{k''\le k'-2}2^{\f {j''}2}2^{k''}\|\Delta_{j''}^v\dot\Delta_{k''}^ha\|_{L^2}\\
&\le& C2^{(1-\sigma_1)k}\|a\|_{H^{\sigma_1,s}}, \eeno from which, it
follows that \ben\label{I-estimate} \|I\|_{L^2}&\le&
C2^{(1-\sigma_1)k}\|a\|_{H^{\sigma_1,s}}\sum_{j'\ge
j-4}\sum_{|k'-k|\le 4}
\|\Delta_{j'}^v\dot\Delta_{k'}^hb\|_{L^2}\nonumber\\
&\le&Cc_{j,k}2^{-js}2^{(1-\sigma_1-\sigma_2)k}\|a\|_{H^{\sigma_1,s}}\|b\|_{H^{\sigma_2,s}}.
\een Similarly, we have \beno II=\sum_{j'\ge j-4}\sum_{|k'-k|\le
4}\Delta_j^v\dot\Delta_k^h
(\Delta_{j'}^vS_{k'-1}^hbS_{j'+2}^v\dot\Delta_{k'}^ha). \eeno Then
we get by Lemma \ref{Lem:Berstein} that \ben\label{II-estimate}
\|II\|_{L^2}
&\le& C\sum_{j'\ge j-4}\sum_{|k'-k|\le 4}\|\Delta_{j'}^vS_{k'-1}^hb\|_{L^2_{x_3}L^\oo_{x_h}}\|S_{j'+2}^v\dot\Delta_{k'}^ha\|_{L^2_{x_h}L^\oo_{x_3}}\nonumber\\
&\le&
C2^{-js}2^{(1-\sigma_1-\sigma_2)k}\|a\|_{H^{\sigma_1,s}}\|b\|_{H^{\sigma_2,s}}
\sum_{j'\ge j-4}\sum_{|k'-k|\le 4}2^{-(j'-j)s}c_{k'}c_{j'}\nonumber\\
&\le&
Cc_{j,k}2^{-js}2^{(1-\sigma_1-\sigma_2)k}\|a\|_{H^{\sigma_1,s}}\|b\|_{H^{\sigma_2,s}}.
\een Now, let us turn to $III$. We have \beno III=\sum_{j'\ge
j-4}\sum_{k',k''\ge k-2;|k'-k''|\le 1}\Delta_j^v\dot\Delta_k^h
(S_{j'+2}^v\dot\Delta_{k'}^ha\Delta_{j'}^v\dot\Delta_{k''}^hb).
\eeno So, we have by Lemma \ref{Lem:Berstein} that
\ben\label{III-estimate}
&&\|III\|_{L^2}\nonumber\\
&&\le C\sum_{j'\ge j-4}\sum_{k',k''\ge k-2;|k'-k''|\le 1}2^k\|S_{j'+2}^v\dot\Delta_{k'}^ha\Delta_{j'}^v\dot\Delta_{k''}^hb\|_{L^2_{x_3}L^1_{x_h}}\nonumber\\
&&\le C\sum_{j'\ge j-4}\sum_{k',k''\ge k-2;|k'-k''|\le
1}2^k\|S_{j'+2}^v\dot\Delta_{k'}^ha\|_{L^\infty_{x_3}L^2_{x_h}}
\|\Delta_{j'}^v\dot\Delta_{k''}^hb\|_{L^2}\nonumber\\
&&\le C2^{-js}2^{(1-\sigma_1-\sigma_2)k}\|a\|_{H^{\sigma_1,s}}\|b\|_{H^{\sigma_2,s}}\sum_{j'\ge j-4}\sum_{k'\ge k-2}2^{-(\sigma_1+\sigma_2)(k'-k)}2^{-(j'-j)s}c_{k'}c_{j'}\nonumber\\
&&\le
Cc_{j,k}2^{-js}2^{(1-\sigma_1-\sigma_2)k}\|a\|_{H^{\sigma_1,s}}\|b\|_{H^{\sigma_2,s}}.
\een

Summing up (\ref{I-estimate})-(\ref{III-estimate}), we obtain \beno
\|\Delta_j^v\dot\Delta_k^h(R_ab)\|_{L^2}\le
Cc_{j,k}2^{-js}2^{(1-\sigma_1-\sigma_2)k}\|a\|_{H^{\sigma_1,s}}\|b\|_{H^{\sigma_2,s}}.
\eeno The lemma is proved in the case when the function $\Psi$ is
identically 0. In order to treat the general case, we only need to
notice the fact that \beno
|\cF\bigl[\Delta_j\dot\Delta_k^h(R_ab)\bigr]_\Psi(\xi)|\le
\cF\bigl[\Delta_j\dot\Delta_k^h(R_{a^+_\Psi}{b^+_\Psi})\bigr](\xi).
\eeno

This finishes the proof of Lemma
\ref{Lem:localproduct}.\ef\vspace{0.1cm}

As a consequence of Lemma \ref{Lem:localproduct} and
(\ref{equ:equinorm}), we have

\bthm{Lemma}\label{Lem:product} Let $s>\f12,\sigma_1,\sigma_2<1$ and
$\sigma_1+\sigma_2>0$. Assume that $a_\Psi\in H^{\sigma_1,s}$ and
$b_\Psi\in H^{\sigma_2,s}$. Then there holds
\beno
\|(ab)_\Psi\|_{H^{\sigma_1+\sigma_2-1,s}}\le C\|a_\Psi\|_{H^{\sigma_1,s}}\|b_\Psi\|_{H^{\sigma_2,s}}.
\eeno

\ethm

\section{Classical analytical-type estimates}

In this section, we prove Proposition \ref{prop:tht}. In this part, we don't need to use any regularizing effect
from the analyticity, but only the fact that the $e^{\Phi(t,\xi_3)}$ is a sublinear function.

Notice that $\p_t v_{\Phi}+\lambda\dot \tht(t)|D_3|v_\Phi=(\p_tv)_\Phi$, we find from (\ref{equ:RNS}) that
\begin{equation}\label{equ:RNS-Phi}
\left\{
\begin{array}{ll}
\p_tv^h_\Phi+\lambda\dot \tht(t)|D_3|v^h_\Phi-\Delta_h v^h_\Phi-\veps^2 \p_3^2 v^h_\Phi+\veps^\f12 (v\cdot \na v^h)_\Phi=-\na_h q_\Phi,\\
\p_tv^3_\Phi+\lambda\dot \tht(t)|D_3|v^3_\Phi-\Delta_h v^3_\Phi-\veps^2 \p_3^2 v^3_\Phi+\veps^\f12(v\cdot\na v^3)_\Phi=-\veps^2\p_3 q_\Phi,\\
\textrm{div\,}v_\Phi= 0,\\
v_\Phi(0)=e^{a|D_3|}v_{0}(x).
\end{array}
\right.
\end{equation}

{\bf Step 1. Estimates on the vertical component $v^3_\Phi$}\vspace{0.2cm}

Note that $\dot \tht(t)\ge 0$, we get from the second equation of (\ref{equ:RNS-Phi}) that
\ben
&&\f 12 \f d{dt}\|v_\Phi^3(t)\|_{H^{-\f12,s}}^2
+\|\na_h v_\Phi^3(t)\|_{H^{-\f12,s}}^2+\|\veps \p_3v_\Phi^3(t)\|_{H^{-\f12,s}}^2\nonumber\\
&&\le -\veps^\f12\bigl((v^h\cdot\na_h v^3)_\Phi, v^3_\Phi\bigr)_{H^{-\f12,s}}
+\veps^\f12\bigl((v^3\textrm{div}_hv^h)_\Phi, v^3_\Phi\bigr)_{H^{-\f12,s}}-\veps^2\bigl(\p_3 q_\Phi,v^3_\Phi\bigr)_{H^{-\f12,s}}\nonumber\\
&&\eqdefa I+II+III.\label{equ:venergy-diff}
\een
Here we used the fact that  $\textrm{div\,}v=0$ such that
\beno
v\cdot\na v^3=v^h\cdot\na_h v^3-v^3\textrm{div}_hv^h.
\eeno
For $II$, Lemma \ref{Lem:product} applied gives
\ben\label{equ:venergy-II}
|II|&\le& \veps^{\f12} \|(v^3\textrm{div}_hv^h)_\Phi\|_{H^{-\f12,s}}\|v^3_\Phi\|_{H^{-\f12,s}}\nonumber\\
&\le& C \veps^{\f12}\|v^3_\Phi\|_{H^{\f12,s}}\|\nabla_h v^h_\Phi\|_{H^{0,s}}\|v^3_\Phi\|_{H^{-\f12,s}}\nonumber\\
&\le&  \f1 {100}\|v^3_\Phi\|^2_{H^{\f12,s}}+C\veps\|\nabla_hv^h_\Phi\|^2_{H^{0,s}}\|v^3_\Phi\|_{H^{-\f12,s}}^2.
\een
For $I$, we get by integration by parts that
\beno
I=\veps^\f12\bigl((\textrm{div}_hv^hv^3)_\Phi, v^3_\Phi\bigr)_{H^{-\f12,s}}
+\veps^\f12\bigl((v^h v^3)_\Phi, \na_hv^3_\Phi\bigr)_{H^{-\f12,s}}\eqdefa I_1+I_2.
\eeno
As in (\ref{equ:venergy-II}), we have
\ben\label{equ:venergy-I1}
|I_1|\le  \f1 {100}\|v^3_\Phi\|^2_{H^{\f12,s}}+C\veps\|\nabla_hv^h_\Phi\|^2_{H^{0,s}}\|v^3_\Phi\|_{H^{-\f12,s}}^2,
\een
and by Lemma \ref{Lem:product},
\ben\label{equ:venergy-I2}
|I_2|&\le& \veps^{\f12}\|(v^3v^h)_\Phi\|_{H^{-\f12,s}}\|\na_hv^3_\Phi\|_{H^{-\f12,s}}\nonumber\\
&\le& C \veps^{\f12}\|v^h_\Phi\|_{H^{0,s}}\|v^3_\Phi\|_{H^{\f12,s}}\|\na_hv^3_\Phi\|_{H^{-\f12,s}}\nonumber\\
&\le& C\veps\|v^h_\Phi\|_{H^{0,s}}^2\|v^3_\Phi\|_{H^{\f12,s}}^2+\f 1 {100}\|\na_hv^3_\Phi\|_{H^{-\f12,s}}^2.
\een

Now, we turn to the estimates of the pressure. Recall that the pressure verifies
$$-\Delta_\veps p=\veps^{\f12}\big[\partial_i \partial_j (v^i v^j) +\partial_i \partial_3( v^i v^3)-2
\partial_3(v^3 \text{div\,}_h v^h)\big].$$
Here and in what follows the index $i,j$ run from 1 to 2. Thus, we can write $p=p^1+p^2+p^3$ with
\begin{equation} \label{equ:pressure-decom}
\begin{split}
&p^1=\veps^{\f12}(-\Delta_\veps)^{-1}\partial_i \partial_j(v^i v^j),\\
&p^2=\veps^{\f12}(-\Delta_\veps)^{-1}\partial_i \partial_3(v^i v^3),\\
&p^3=-2\veps^{\f12}(-\Delta_\veps)^{-1}\partial_3(v^3\text{div\,}_h v^h).
 \end{split}
 \end{equation}
We get by integration by parts that
$$
\veps^2(\partial_3 p^1_\Phi,
v^3_\Phi)_{H^{-\f12,s}}=-\veps(p^1_\Phi, \veps\partial_3
v^3_\Phi)_{H^{-\f12,s}}\leq C\veps^2
\|p^1_\Phi\|_{H^{-\f12,s}}^2+\f {1} {100}\|\veps\partial_3
v^3_\Phi\|_{H^{-\f12,s}}^2,
$$
which together with the fact that the operator $\partial_i\partial_j(-\Delta_\veps)^{-1}$
is bounded on $H^{\sigma, s}$ and Lemma \ref{Lem:product} implies that
\ben\label{equ:venergy-p1}
\veps^2(\partial_3 p^1_\Phi, v^3_\Phi)_{H^{-\f12,s}}&\leq& C\veps^{3} \|(v^h\otimes v^h)_\Phi\|_{H^{-\f12,s}}^2
+\f 1 {100}\|\veps\partial_3 v^3_\Phi\|_{H^{-\f12,s}}^2\nonumber\\
&\leq& C\veps^2\|\veps^\f12v^h_\Phi\|_{H^{\f12,s}}^2\|v^h_\Phi\|_{H^{0,s}}^2+
\f 1 {100}\|\veps\partial_3 v^3_\Phi\|_{H^{-\f12,s}}^2.
\een
For the term containing $p_2$, we get by integration by parts that
\beno
\veps^2(\partial_3 p^2_\Phi, v^3_\Phi)_{H^{-\f12,s}}=-\veps^\f12(\veps^2\partial_3^2 (-\Delta_\veps)^{-1}(v^i v^3)_\Phi, \p_iv^3_\Phi)_{H^{-\f12,s}},
\eeno
then using the fact that $(\veps\partial_3)^2(-\Delta_\veps)^{-1}$
is bounded  on $H^{\sigma, s}$ and Lemma \ref{Lem:product}, we have
\ben\label{equ:venergy-p2}
\veps^2(\partial_3 p^2_\Phi, v^3_\Phi)_{H^{-\f12,s}}&\leq& C\veps^{\f12}\|(v^3 v^h)_\Phi\|_{H^{-\f12,s}}\|\nabla_hv^3_\Phi\|_{H^{-\f12,s}}\nonumber \\
&\leq& C\veps\|v^3_\Phi\|_{H^{\f12,s}}^2\|v^h_\Phi\|_{H^{0,s}}^2+\f 1{100}\|\nabla_h
v^3_\Phi\|_{H^{-\f12,s}}^2.
\een
For the last term coming from $p_3$, we use again the fact that
$(\veps\partial_3)^2(-\Delta_\veps)^{-1}$ is  bounded on $H^{\sigma, s}$
and obtain
\ben\label{equ:venergy-p3}
\veps^2(\partial_3 p^3_\Phi, v^3_\Phi)_{H^{-\f12,s}}&\leq& C\veps^{\f12}\|(v^3\text{div\,}v^h)_\Phi\|_{H^{-\f12,s}}\|v^3_\Phi\|_{H^{-\f1 2,s}}\nonumber\\
&\leq& C\veps^{\f12}\|v^3_\Phi\|_{H^{\f12,s}}\|\nabla_h v^h_\Phi\|_{H^{0,s}}\|v^3_\Phi\|_{H^{-\f12,s}}\nonumber\\
&\leq& C\veps \|\nabla_hv^h_\Phi\|^2_{H^{0,s}}\|v^3_\Phi\|_{H^{-\f12,s}}^2+\f 1{100}\|v^3_\Phi\|_{H^{\f12,s}}^2.
\een

Summing up (\ref{equ:venergy-diff})-(\ref{equ:venergy-I2}) and (\ref{equ:venergy-p1})-(\ref{equ:venergy-p3}), we obtain
\ben\label{equ:venergy}
&&\f d{dt}\|v_\Phi^3(t)\|_{H^{-\f12,s}}^2
+\|v_\Phi^3(t)\|_{H^{\f12,s}}^2\nonumber\\&&\le C\|\nabla_hv^h_\Phi\|^2_{H^{0,s}}\|v^3_\Phi\|_{H^{-\f12,s}}^2+
C(\|v^3_\Phi\|_{H^{\f12,s}}^2+\|\veps^\f12v^h_\Phi\|_{H^{\f12,s}}^2)\|v^h_\Phi\|_{H^{0,s}}^2.
\een
Here we used the fact that
\beno
\|\na_h v_\Phi^3\|_{H^{-\f12,s}}^2\thickapprox\|v_\Phi^3\|_{H^{\f12,s}}^2.
\eeno

{\bf Step 2. Estimates on the horizontal component $v^h_\Phi$}\vspace{0.2cm}

From the first equation of (\ref{equ:RNS-Phi}), we infer that
\ben\label{equ:henergy-diff}
&&\f 12 \f d{dt}\|\veps^\f12v_\Phi^h(t)\|_{H^{-\f12,s}}^2
+\|\veps^\f12\na_h v_\Phi^h(t)\|_{H^{-\f12,s}}^2+\veps\|\veps \p_3v_\Phi^h(t)\|_{H^{-\f12,s}}^2\nonumber\\
&&\quad\le -\veps\bigl((v\cdot\na v^h)_\Phi,
\veps^\f12v^h_\Phi\bigr)_{H^{-\f12,s}} -\veps\bigl(\na_h
q_\Phi,v^h_\Phi\bigr)_{H^{-\f12,s}}\eqdefa I+II.
\een
We rewrite $I$ as \beno I=-\veps\bigl((v^h\cdot\na_h v^h)_\Phi,
\veps^\f12v^h_\Phi\bigr)_{H^{-\f12,s}}
-\veps\bigl((v^3\p_3v^h)_\Phi,
\veps^\f12v^h_\Phi\bigr)_{H^{-\f12,s}}\eqdefa I_1+I_2. \eeno
Lemma \ref{Lem:product} applied gives
\ben\label{equ:henergy-I1}
|I_1|&\le& \veps\|(v^h\na_hv^h)_\Phi\|_{H^{-\f12,s}}\|\veps^\f12v^h_\Phi\|_{H^{-\f12,s}}\nonumber\\
&\le& C \veps\|v^h_\Phi\|_{H^{\f12,s}}\|\nabla_h v^h_\Phi\|_{H^{0,s}}\|\veps^\f12v^h_\Phi\|_{H^{-\f12,s}}\nonumber\\
&\le&  \f1 {100}\|\veps^\f12v^h_\Phi\|^2_{H^{\f12,s}}+C\veps\|\nabla_hv^h_\Phi\|^2_{H^{0,s}}\|\veps^\f12v^h_\Phi\|_{H^{-\f12,s}}^2.
\een
For $I_2$, we use integration by parts and $\text{div\,}v=0$ to get
\beno
I_2&=&-\veps\bigl((\textrm{div}_hv^hv^h)_\Phi, \veps^\f12v^h_\Phi\bigr)_{H^{-\f12,s}}
+\bigl((v^h v^3)_\Phi, \veps^\f12\veps\p_3v^h_\Phi\bigr)_{H^{-\f12,s}},\\
&\eqdefa& I_{21}+I_{22}.
\eeno
As in (\ref{equ:henergy-I1}), we have
\ben\label{equ:henergy-I21}
|I_{21}|\le  \f1 {100}\|\veps^\f12v^h_\Phi\|^2_{H^{\f12,s}}+C\veps\|\nabla_hv^h_\Phi\|^2_{H^{0,s}}\|\veps^\f12v^h_\Phi\|_{H^{-\f12,s}}^2,
\een
and by Lemma \ref{Lem:product},
\ben\label{equ:henergy-I22}
|I_{22}|&\le& \|(v^3v^h)_\Phi\|_{H^{-\f12,s}}\veps^\f12\|\veps\p_3v^h_\Phi\|_{H^{-\f12,s}}\nonumber\\
&\le& C\|v^h_\Phi\|_{H^{0,s}}\|v^3_\Phi\|_{H^{\f12,s}}\veps^\f12\|\veps\p_3v^h_\Phi\|_{H^{-\f12,s}}\nonumber\\
&\le& \f1 {100}\veps\|\veps\p_3v^h_\Phi\|^2_{H^{\f12,s}}+C\|v^h_\Phi\|^2_{H^{0,s}}\|v^3_\Phi\|_{H^{\f12,s}}^2.
\een

In order to deal with the pressure, we write $p=p^1+p^2+p^3$ with $p^1,p^2,p^3$ defined by (\ref{equ:pressure-decom}).
Using the fact that the operator $\partial_i\partial_j(-\Delta_\veps)^{-1}$
is bounded on $H^{\sigma, s}$ and Lemma \ref{Lem:product}, we have
\ben\label{equ:henergy-p1}
\veps(\na_h p^1_\Phi, v^h_\Phi)_{H^{-\f12,s}}
&=& -\veps((-\Delta_\veps)^{-1}\partial_i \partial_j(v^i v^j)_\Phi,\veps^\f12\text{div}_hv^h_\Phi)_{H^{-\f12,s}}\nonumber\\
&\leq& C\veps\|(v^h\otimes v^h)_\Phi\|_{H^{-\f12,s}}\|\veps^\f12\na_hv^h_\Phi\|_{H^{-\f12,s}}\nonumber\\
&\leq& C\veps^\f12\|\veps^\f12v^h_\Phi\|_{H^{\f12,s}}\|v^h_\Phi\|_{H^{0,s}}\|\veps^\f12\na_hv^h_\Phi\|_{H^{-\f12,s}}\nonumber\\
&\le& \f 1 {100}\|\veps^\f12\na_hv^h_\Phi\|_{H^{-\f12,s}}^2+C\veps\|\veps^\f12v^h_\Phi\|_{H^{\f12,s}}^2\|v^h_\Phi\|_{H^{0,s}}^2.
\een
For the term coming from $p_2$, we integrate by parts to get
\beno
\veps(\na_h p^2_\Phi, v^h_\Phi)_{H^{-\f12,s}}=-(\veps\p_i\partial_3(-\Delta_\veps)^{-1}(v^i v^3)_\Phi, \veps^\f12\text{div\,}_hv^h_\Phi)_{H^{-\f12,s}},
\eeno
then note that $\veps\partial_3\p_i(-\Delta_\veps)^{-1}$ is  bounded  on $H^{\sigma, s}$,  we get by Lemma \ref{Lem:product} that
\ben\label{equ:henergy-p2}
\veps(\na_h p^2_\Phi, v^h_\Phi)_{H^{-\f12,s}}&\leq& C\|(v^3 v^h)_\Phi\|_{H^{-\f12,s}}\|\veps^\f12\nabla_hv^h_\Phi\|_{H^{-\f12,s}}\nonumber \\
&\le& C\|v^3_\Phi\|_{H^{\f12,s}}\|v^h_\Phi\|_{H^{0,s}}\|\veps^\f12\nabla_h
v^h_\Phi\|_{H^{-\f12,s}}\nonumber\\
&\leq& C\|v^3_\Phi\|_{H^{\f12,s}}^2\|v^h_\Phi\|_{H^{0,s}}^2+\f 1{100}\|\veps^\f12\nabla_h
v^h_\Phi\|_{H^{-\f12,s}}^2.
\een
Similarly, we have
\ben\label{equ:henergy-p3}
\veps(\na_h p^3_\Phi, v^h_\Phi)_{H^{-\f12,s}}&\leq& C\|(v^3\text{div\,}_hv^h)_\Phi\|_{H^{-\f12,s}}\|\veps^\f12v^h_\Phi\|_{H^{-\f1 2,s}}\nonumber\\
&\leq& C\|v^3_\Phi\|_{H^{\f12,s}}\|\nabla_h v^h_\Phi\|_{H^{0,s}}\|\veps^\f12v^h_\Phi\|_{H^{-\f12,s}}\nonumber\\
&\leq&\f 1{100}\|v^3_\Phi\|_{H^{\f12,s}}^2+ C\|\nabla_hv^h_\Phi\|^2_{H^{0,s}}\|\veps^\f12v^h_\Phi\|_{H^{-\f12,s}}^2.
\een

Summing up (\ref{equ:henergy-diff})-(\ref{equ:henergy-p3}), we obtain
\ben\label{equ:henergy}
&&\f d{dt}\|\veps^\f12v_\Phi^h(t)\|_{H^{-\f12,s}}^2+\|\veps^\f12 v_\Phi^h(t)\|_{H^{\f12,s}}^2-\f1{50}\|v^3_\Phi\|_{H^{\f12,s}}^2\nonumber\\
&&\le
C\|\nabla_hv^h_\Phi\|^2_{H^{0,s}}\|\veps^\f12v^h_\Phi\|_{H^{-\f12,s}}^2+
C\bigl(\|v^3_\Phi\|_{H^{\f12,s}}^2+\|\veps^\f12v^h_\Phi\|_{H^{\f12,s}}^2\bigr)\|v^h_\Phi\|_{H^{0,s}}^2.
\een

{\bf Step 3. Estimate on the function $\tht(t)$}\vspace{0.2cm}

Combining (\ref{equ:venergy}) with (\ref{equ:henergy}), we obtain
\beno
&&\f d{dt}\bigl(\|\veps^\f12v_\Phi^h(t)\|_{H^{-\f12,s}}^2+\|v_\Phi^3(t)\|_{H^{-\f12,s}}^2\bigr)+
\bigl(\|\veps^\f12 v_\Phi^h(t)\|_{H^{\f12,s}}^2+\|v_\Phi^3(t)\|_{H^{\f12,s}}^2\bigr)\nonumber\\
&&\quad\le C\|\nabla_hv^h_\Phi\|^2_{H^{0,s}}\bigl(\|\veps^\f12v^h_\Phi\|_{H^{-\f12,s}}^2+\|v^3_\Phi\|_{H^{-\f12,s}}^2\bigr)\nonumber\\
&&\qquad+C\bigl(\|\veps^\f12v^h_\Phi\|_{H^{\f12,s}}^2+\|v^3_\Phi\|_{H^{\f12,s}}^2\bigr)\|v^h_\Phi\|_{H^{0,s}}^2,
\eeno
from which and Gronwall's inequality, it follows that
\beno
&&\|\veps^\f12v_\Phi^h(t)\|_{H^{-\f12,s}}^2+\|v_\Phi^3(t)\|_{H^{-\f12,s}}^2+
\int_0^t\bigl(\|\veps^\f12 v_\Phi^h(\tau)\|_{H^{\f12,s}}^2+\|v_\Phi^3(\tau)\|_{H^{\f12,s}}^2\bigr)d\tau\nonumber\\
&&\quad\le \exp\bigl(C\int_0^t\|\nabla_hv^h_\Phi(\tau)\|^2_{H^{0,s}}d\tau\bigr)\Bigl[\|e^{a|D_3|}v_0\|_{H^{-\f12,s}}^2+\nonumber\\
&&\qquad\quad+C\int_0^t\bigl(\|\veps^\f12v^h_\Phi(\tau)\|_{H^{\f12,s}}^2+\|v^3_\Phi(\tau)\|_{H^{\f12,s}}^2\bigr)\|v^h_\Phi(\tau)\|_{H^{0,s}}^2d\tau\Bigr].
\eeno
In particular, we have
\beno
\tht(t)\le\exp\bigl(C\Psi(t)\bigr)\Bigl[\|e^{a|D_3|}v_0\|_{H^{-\f12,s}}^2+
C\int_0^t\dot\tht(\tau)\Psi(\tau)d\tau\Bigr].
\eeno

This finishes the proof of Proposition \ref{prop:tht}.\ef

\section{Regularizing effect du analyticity}

Let's now prove Proposition \ref{prop:Psi}. Here we will encounter two kinds of bad
terms, where we lose a vertical derivative. The first one is
$(v^3\partial_3 v^h)_\Phi$. We will see that in an energy estimate, we
have no loss of vertical derivative in this term (by integrating by
parts, using commutators and of course $\partial_3 v^3=-\text{div\,}_h
v^h$). In the term $\nabla p$, we really lose a vertical derivative.

{\bf Step 1. Estimates on the horizontal component $v^h_\Phi$}\vspace{0.2cm}

Let us recall that $v_\Phi^h$ verifies the equations
\beno
\p_tv^h_\Phi+\lambda\dot \tht(t)|D_3|v^h_\Phi-\Delta_h v^h_\Phi-\veps^2 \p_3^2 v^h_\Phi+\veps^\f12 (v\cdot \na v^h)_\Phi=-\na_h q_\Phi.
\eeno
Note that $\dot\theta \geq 0$, we perform an energy estimate
in $H^{0,s}$ to obtain
\ben\label{equ:Psi-diff}
&&\f12\frac{d}{dt}\|v^h_\Phi\|^2_{H^{0,s}}+\lambda\dot\theta(t)\|v^h_\Phi\|^2_{H^{0,s+1/2}}+\|\na_hv^h_\Phi\|_{H^{0,s}}^2
+\|\veps\p_3v^h_\Phi\|_{H^{0,s}}^2\nonumber\\
&&\leq \veps^\f12((v^h\otimes v^h)_\Phi, \na_h v^h_\Phi)_{H^{0,s}}-\veps^\f12(\p_3(v^3v^h)_\Phi,v^h_\Phi)_{H^{0,s}}-(\nabla_h p_\Phi,v^h_\Phi)_{H^{0,s}}\nonumber\\
&&\eqdefa I+II+III.
\een
We get by Lemma \ref{Lem:product} and the interpolation that
\ben\label{equ:Psi-I}
|I|&\leq& C\veps^\f12\|(v_h\otimes v^h)_\Phi\|_{H^{0,s}}\|\na_h v^h_\Phi\|_{H^{0,s}}\nonumber\\
&\leq& C\veps^\f12\|v^h_\Phi\|_{H^{\f12,s}}\|v^h_\Phi\|_{H^{\f12,s}}\|\nabla_h v^h_\Phi\|_{H^{0,s}}\nonumber\\
&\le& C\veps^\f12\|v^h_\Phi\|_{H^{0,s}}\|v^h_\Phi\|_{H^{1,s}}\|\nabla_h v^h_\Phi\|_{H^{0,s}}\nonumber\\
&\le& C\veps\|v^h_\Phi\|_{H^{0,s}}^2\|v^h_\Phi\|_{H^{1,s}}^2+\f 1 {100}\|\nabla_h v^h_\Phi\|_{H^{0,s}}^2.
\een
In order to estimate $II$, we use Bony's decomposition (\ref{VBony}) to rewrite it as
\beno
II=-\veps^\f12(\p_3(T^v_{v^h}v_3)_\Phi,v^h_\Phi)_{H^{0,s}}-\veps^\f12(\p_3(R^v_{v_3}v^h)_\Phi,v^h_\Phi)_{H^{0,s}}
\eqdefa II_1+II_2.
\eeno
From the proof of  Lemma \ref{Lem:localproduct}, it is easy to find that
\ben\label{equ:Psi-II2}
|II_2|&\le& C\||D_3|^\f12(R^v_{v^3}v^h)_\Phi\|_{H^{-\f12,s}}\||\veps D_3|^\f12|D_h|^\f12v^h_\Phi\|_{H^{0,s}}\nonumber\\
&\le& C\|v^3_{\Phi}\|_{H^{\f12,s}}\|v^h_\Phi\|_{H^{0,s+\f12}}\|\na_\veps v^h\|_{H^{0,s}}\nonumber\\
&\le& C\|v^3_{\Phi}\|_{H^{\f12,s}}^2\|v^h_\Phi\|_{H^{0,s+\f12}}^2+\f1 {100}\|\na_\veps v^h\|_{H^{0,s}}^2.
\een
Due to $\text{div\,}v=0$, we rewrite $II_1$ as
\beno
II_1=\veps^\f12((T^v_{v^h}\text{div\,}v^h)_\Phi,v^h_\Phi)_{H^{0,s}}-\veps^\f12((T^v_{\p_3v^h}v^3)_\Phi,v^h_\Phi)_{H^{0,s}}
\eqdefa II_{11}+II_{12}.
\eeno
Using Lemma \ref{Lem:localproduct}, we have
\ben\label{equ:Psi-II11}
|II_{11}|&\le& \veps^\f12\|(T^v_{v^h}\text{div\,}v^h)_\Phi\|_{H^{-\f12,s}}\|v^h_\Phi\|_{H^{\f12,s}}\nonumber\\
&\le& C\veps^\f12\|v^h_\Phi\|_{H^{\f12,s}}\|\nabla_h v^h_\Phi\|_{H^{0,s}}\|v^h_\Phi\|_{H^{\f12,s}}\nonumber\\
&\le& C\veps\|v^h_\Phi\|_{H^{0,s}}^2\|v^h_\Phi\|_{H^{1,s}}^2+\f 1 {100}\|\nabla_h v^h_\Phi\|_{H^{0,s}}^2.
\een
From the proof of  Lemma \ref{Lem:localproduct}, we can conclude that
\ben\label{equ:Psi-II12}
|II_{12}| &\le& C\|v^3_{\Phi}\|_{H^{\f12,s}}\|v^h_\Phi\|_{H^{0,s+\f12}}\|\na_\veps v^h\|_{H^{0,s}}\nonumber\\
&\le& C\|v^3_{\Phi}\|_{H^{\f12,s}}^2\|v^h_\Phi\|_{H^{0,s+\f12}}^2+\f1 {100}\|\na_\veps v^h\|_{H^{0,s}}^2.
\een

We next turn to the estimate of the pressure. Recall that $p=p^1+p^2+p^3$ with $p^1,p^2,p^3$ defined by (\ref{equ:pressure-decom}).
Using the fact that $(-\Delta_\veps)^{-1}\partial_i\partial_j$
is bounded on $H^{\sigma, s}$ and Lemma \ref{Lem:product}, we get
\ben\label{equ:Psi-p1}
(\na_h p^1_\Phi, v^h_\Phi)_{H^{0,s}}
&=& -\veps^\f12((-\Delta_\veps)^{-1}\partial_i \partial_j(v^i v^j)_\Phi,\text{div\,}v^h_\Phi)_{H^{0,s}}\nonumber\\
&\leq& C\veps^\f12\|(v^h\otimes v^h)_\Phi\|_{H^{0,s}}\|\na_hv^h_\Phi\|_{H^{0,s}}\nonumber\\
&\leq& C\veps^\f12\|v^h_\Phi\|_{H^{\f12,s}}\|v^h_\Phi\|_{H^{\f12,s}}\|\na_hv^h_\Phi\|_{H^{0,s}}\nonumber\\
&\le&  C\veps\|v^h_\Phi\|_{H^{0,s}}^2\|v^h_\Phi\|_{H^{1,s}}^2+\f 1 {100}\|\nabla_h v^h_\Phi\|_{H^{0,s}}^2.
\een
Notice that $\partial_i\partial_j(-\Delta_\veps)^{-1}$
is bounded on $H^{\sigma, s}$, then exactly as in the estimate of $II$, we can obtain
\ben\label{equ:Psi-p2}
&&(\na_h p^2_\Phi, v^h_\Phi)_{H^{0,s}}\nonumber\\&&\leq C\|v^3_{\Phi}\|_{H^{\f12,s}}^2\|v^h_\Phi\|_{H^{0,s+\f12}}^2+
C\veps\|v^h_\Phi\|_{H^{0,s}}^2\|v^h_\Phi\|_{H^{1,s}}^2+\f1 {100}\|\na_\veps v^h\|_{H^{0,s}}^2.
\een
We write
\beno
\nabla_h p_3= -2\partial_3|D_3|^{-\f12}\bigl(\na_h|D_h|^{\f12}|\veps D_3|^{1/2}(-\Delta_\epsilon)^{-1}\bigr)|D_h|^{-\f12}(v^3\text{div\,}_h v^h)
\eeno
thus,
\beno
(\na_h p^3_\Phi, v^h_\Phi)_{H^{0,s}}=-2(\bigl(\na_h|D_h|^{\f12}|\veps D_3|^{\f12}(-\Delta_\epsilon)^{-1}\bigr)|D_h|^{-\f12}(v^3\text{div\,}_h v^h),
\partial_3|D_3|^{-\f12}v^h_\Phi)_{H^{0,s}}.
\eeno
Note that $\na_h|D_h|^{\f12}|\veps D_3|^{\f12}(-\Delta_\epsilon)^{-1}$ is a bounded operator on $H^{\sigma,s}$,
we get by Lemma \ref{Lem:product} that
\ben\label{equ:Psi-p3}
(\nabla_h p^3_\Phi, v^h_\Phi)_{H^{0,s}}&\le& C\||D_h|^{-1/2}(v^3\text{div\,}_h v^h)\|_{H^{0,s}}\|\partial_3|D_3|^{-\f12}v^h_\Phi\|_{H^{0,s}}\nonumber\\
&\leq& C\|v^3_\Phi\|_{H^{\f12,s}}\|\nabla_h v^h_\Phi\|_{H^{0,s}}\|v^h\|_{H^{0,s+1/2}}\nonumber\\
&\leq& C\|v^3_\Phi\|_{H^{\f12,s}}^2\|v^h_\Phi\|_{H^{0,s+1/2}}^2+\f 1{100}\|\nabla_h v^h_\Phi\|^2_{H^{0,s}}.
\een

Summing up (\ref{equ:Psi-diff})-(\ref{equ:Psi-p3}), we get by taking $\lambda$ big enough that
\ben\label{equ:hPsi-energy}
\frac{d}{dt}\|v^h_\Phi(t)\|^2_{H^{0,s}}+\|\na_hv^h_\Phi(t)\|_{H^{0,s}}^2
\le C\|v^h_\Phi\|_{H^{0,s}}^2\|\na_hv^h_\Phi\|_{H^{0,s}}^2.
\een

{\bf Step 2. Estimates on the vertical component $v^3_\Phi$}\vspace{0.2cm}

Recall that $v_\Phi^3$ verifies the equation
\beno
\p_tv^3_\Phi+\lambda\dot \tht(t)|D_3|v^3_\Phi-\Delta_h v^3_\Phi-\veps^2 \p_3^2 v^3_\Phi+\veps^\f12 (v\cdot \na v^3)_\Phi
=-\veps^2\p_3 q_\Phi.
\eeno
We perform an energy estimate in $H^{0,s}$ to obtain
\ben\label{equ:vPsi-diff}
&&\f12\frac{d}{dt}\|v^3_\Phi\|^2_{H^{0,s}}+\|\na_hv^3_\Phi\|_{H^{0,s}}^2
+\|\veps\p_3v^3_\Phi\|_{H^{0,s}}^2\nonumber\\
&&\leq -\veps^\f12((v^h\cdot\na_h v^3)_\Phi, v^3_\Phi)_{H^{0,s}}+\veps^\f12((v^3\text{div\,}_hv^h)_\Phi,v^3_\Phi)_{H^{0,s}}
-\veps^2(\p_3 p_\Phi,v^3_\Phi)_{H^{0,s}}\nonumber\\
&&\eqdefa I+II+III.
\een
Using Lemma \ref{Lem:product} and the inerpolation, we have
\ben\label{equ:vPsi-I}
|I|&\leq& C\veps^\f12\|v^h_\Phi\|_{H^{\f12,s}}\|\na_h v^3_\Phi\|_{H^{0,s}}\|v^3_\Phi\|_{H^{\f12,s}}\nonumber\\
&\le& C\veps\bigl(\|v^h_\Phi\|_{H^{0,s}}^2\|v^h_\Phi\|_{H^{1,s}}^2+\|v^3_\Phi\|_{H^{0,s}}^2\|v^3_\Phi\|_{H^{1,s}}^2\bigr)
+\f 1 {100}\|\nabla_h v^3_\Phi\|_{H^{0,s}}^2,
\een
and
\ben\label{equ:vPsi-II}
|II|&\leq& C\veps^\f12\|v^3_\Phi\|_{H^{\f12,s}}\|\na_h v^h_\Phi\|_{H^{0,s}}\|v^3_\Phi\|_{H^{\f12,s}}\nonumber\\
&\le& C\veps\|v^3_\Phi\|_{H^{0,s}}^2\|\na_h v^h_\Phi\|_{H^{0,s}}^2
+\f 1 {100}\|\nabla_h v^3_\Phi\|_{H^{0,s}}^2.
\een
Using the decomposition (\ref{equ:pressure-decom}), we can similarly obtain
\ben\label{equ:vPsi-III}
|III|\le C\veps\|v_\Phi\|_{H^{0,s}}^2\|\na_hv_\Phi\|_{H^{0,s}}^2
+\f 1 {100}\|\nabla_h v^3_\Phi\|_{H^{0,s}}^2.
\een

Summing up (\ref{equ:vPsi-diff})-(\ref{equ:vPsi-III}), we obtain
\ben\label{equ:vPsi-energy}
\frac{d}{dt}\|v^3_\Phi\|^2_{H^{0,s}}+\|\na_hv^3_\Phi\|_{H^{0,s}}^2
\le  C\veps\|v_\Phi\|_{H^{0,s}}^2\|\na_hv_\Phi\|_{H^{0,s}}^2.
\een

Combining (\ref{equ:hPsi-energy}) with (\ref{equ:vPsi-energy}), we get
\beno
\frac{d}{dt}\|v_\Phi\|^2_{H^{0,s}}+\|\na_hv_\Phi\|_{H^{0,s}}^2
\le C\|v_\Phi\|_{H^{0,s}}^2\|\na_hv_\Phi\|_{H^{0,s}}^2,
\eeno
from which and Gronwall's inequality, we infer that
\beno
\|v_\Phi(t)\|^2_{H^{0,s}}+\int_0^t\|\na_hv_\Phi(\tau)\|_{H^{0,s}}^2d\tau\le \|e^{a|D_3|}v_0\|_{H^{0,s}}^2\exp\bigl(C\int_0^t\|\na_h v_\Phi(\tau)\|^2_{H^{0,s}}d\tau\bigr),
\eeno
that is,
\beno
\Psi(t)\le \|e^{a|D_3|}v_0\|_{H^{0,s}}^2\exp\bigl(C\Psi(t)\bigr).
\eeno

This finishes the proof of Proposition \ref{prop:Psi}.\ef

\section*{acknowledgement}
This work was partly done when Zhifei Zhang was visiting Department of Mathematics of Paris-Sud
University as a Postdoctor Fellowship. He would like to thank the hospitality and support of the Department.

\end{document}